\documentclass[11pt]{article}
\usepackage{amsfonts,amssymb,amsmath}

\usepackage[english]{babel}

\newcommand{\halmos}{\rule{1ex}{1.4ex}}

\makeatletter \@addtoreset{equation}{section} \makeatother

\newtheorem{theorem}{Theorem}[section]
\newtheorem{corollary}{Corollary}[section]
\newtheorem{definition}{Definition}[section]

\newtheorem{lemma}{Lemma}[section]
\newtheorem{proposition}{Proposition}[section]

\newtheorem{remark}{Remark}[section]

\newenvironment{proof}{\noindent {\em Proof}.\,\,\,}
{\hspace*{\fill}$\halmos$\medskip}

\newcommand{\beq}{\begin{eqnarray}}
\newcommand{\eeq}{\end{eqnarray}}

\newcommand{\be}{\begin{equation}}
\newcommand{\ee}{\end{equation}}

\newcommand{\bl}{\begin{lemma}}
\newcommand{\el}{\end{lemma}}

\newcommand{\br}{\begin{remark}}
\newcommand{\er}{\end{remark}}

\newcommand{\bt}{\begin{theorem}}
\newcommand{\et}{\end{theorem}}

\newcommand{\bd}{\begin{definition}}
\newcommand{\ed}{\end{definition}}

\newcommand{\bp}{\begin{proposition}}
\newcommand{\ep}{\end{proposition}}

\newcommand{\bc}{\begin{corollary}}
\newcommand{\ec}{\end{corollary}}

\newcommand{\bpr}{\begin{proof}}
\newcommand{\epr}{\end{proof}}

\newcommand{\bi}{\begin{itemize}}
\newcommand{\ei}{\end{itemize}}

\newcommand{\ben}{\begin{enumerate}}
\newcommand{\een}{\end{enumerate}}

\newcommand{\Z}{\mathbb Z}

\newcommand{\N}{\mathbb N}

\newcommand{\E}{\mathbb E}

\newcommand{\pee}{\mathbb P}

\newcommand{\eps}{\ensuremath{\epsilon}}

\def\b{\beta}

\def\phi{\varphi}

\def\L{{\mathbb L}}

\def\bt{\tilde \b}

\def\1{1}

\def\E{{\Bbb E}}
\def\N{{\Bbb N}}

\def\Z{{\Bbb Z}}

\def\one{\hbox{J}\kern-.2em\hbox{I}}

\begin{document}

\title{{\bf Range of a Transient $2d$-Random Walk}}
\author{ Arnaud Le Ny\footnote{Universit\'e de Paris-Sud 91405 Orsay cedex, France. Email: arnaud.leny@math.u-psud.fr}\footnote{This work has been fulfiled within the ANR ``jeunes chercheurs'' program RANDYMECA.}}
\maketitle


\footnotesize
\begin{center}
{\bf Abstract:} 
\end{center}
We study the range of a planar random walk on a randomly oriented lattice, already known to be transient. We prove that the expectation of the range grows linearly, in both the quenched  (for a.e. orientation) and  annealed ("averaged") cases. We also express the rate of growth in terms of the quenched Green function and eventually prove a weak law of large numbers in the (non-Markovian) annealed case.\normalsize
\section{Preliminaries}
We focus on a particular type of {\em random walk in a random environment} (RWRE), where the environment is inherited from orientations of the lattice on which the walker evolves, providing thus two independent sources of randomness: The horizontal orientations of the lattice and this of  the walk performed on it afterwards, once the realization of the orientation has been fixed. We introduce first an {\em horizontally  oriented} square lattice $\mathbb{L}^\eps$. The orientations $\eps=(\eps_y)_{y \in \Z}$  are families of  i.i.d. Rademacher random variables taking values in the product probability space $(E,\mathcal{E},\rho)=\big(\{-1,+1\},\mathcal{P}(\{-1,+1\}, \frac{1}{2} \delta_{-1} + \frac{1}{2} \delta_{+1} \big)^{\otimes \Z}$. A given horizontal level $y$ is then oriented to the right when $\eps_y=+1$, to the left when $\eps_y=-1$, and induces an horizontally oriented version of $\Z^2$ for every realization of the random field $\eps$:
\bd [Horizontally Oriented Lattice $\mathbb{L}^\eps$]
 Let $\eps=(\eps_y)_{y \in \Z} \in \{\pm 1\}^\Z$. The {\em
oriented lattice}
$\mathbb{L}^\eps=(\mathbb{V},\mathbb{A}^\eps)$ is the (random)
directed graph with (deterministic) vertex set $\mathbb{V}=\Z^2$
and (random) edge set $\mathbb{A}^\eps$ defined by the condition that for $u=(u_1,u_2), v=(v_1,v_2) \in \Z^2$, 
$$
(u,v) \in
\mathbb{A}^\eps \;\Longleftrightarrow \;  v_1=u_1 \; {\rm and} \; v_2=u_2 \pm 1, \; 
{\rm or}
\; v_2=u_2 \;  {\rm and} \; v_1=u_1+ \eps_{u_2}.
$$
\ed
One performs then a {\em simple random walk} (SRW) $M=(M_n)_{n \in \N}$ on $\L^\eps$.  For a given  $\eps$, it is a $\mathbb{Z}^2$-valued Markov chain defined on a probability space $\big(\Omega_0,
 \mathcal{B}_0,\mathbb{P}^{(\eps)}\big)$, starting at the origin $(0,0)$, whose ($\eps$-dependent) transition probabilities are
 defined for all $(u,v) \in \mathbb{V}\times \mathbb{V}$ by

$$
\pee^{(\eps)}[M_{n+1}=v  | M_n=u]=\frac{1}{3}  \; \rm{if} \;  (u,v)   \in \mathbb{A}^\eps, \; 0 \; \; \rm{otherwise.}
$$

An interesting feature is that this SRW has been proven to be {\em transient for almost every orientation} $\eps$ \cite{CP2}. This almost sure approach is referred as the {\em quenched case} and we focus here on a more collective {\em annealed} approach: We consider the law of the process under the joint measure $\pee:= \rho \otimes \pee^{(\eps)}$. Thus, we study  the behavior of the SRW as a discrete-time process on
$$
(\Omega,\mathcal{B},\pee):=\big(E \times \Omega_0, \mathcal{E} \otimes \mathcal{B}_0, \rho \otimes \pee^{(\eps)} \big).
$$
with its  {\em annealed law} $\pee$ formally defined as $\pee= \int_E \pee^{(\eps)} d \rho(\eps).$ We write $\E$ (or $\E^{(\eps)}$ or $\E_\rho$) for the expectation under  $\pee$ (or $\pee^{(\eps)}$ or $\rho$). Due to the non-local character of the orientations, the main drawback of this annealed model is that the walk is {\em not Markovian anymore}. Nevertheless 
\bp 
Under the {\em  annealed} law $\pee$, the process $M$ is {\em reversible}.
\ep
Indeed, consider a trajectory $\omega=(\omega_0, \dots,\omega_n)$ and change $\eps$ into $-\eps$ : It has the law of the reverse trajectory $\omega^*=(\omega_n,\dots,\omega_0)$ and one concludes using the symmetry of the law $\rho$ of $\eps$.\\

Under this annealed law, a non-standard functional limit theorem has been proven in \cite{GPLN2}, while we shall start our study of the range of the random walk thanks to the following estimation of the probability of return to the origin due to Castell {\em et al.} \cite{CGPS}:

\begin{theorem}[Local Limit Theorem \cite{CGPS}] 
There exists a constant $C >0$ such that 
\be \label{LLT}
u_{\rm{cp}}(n):=\pee[M_n =(0,0)] = C \cdot n^{-5/4}  + \circ(n^{-5/4}) \;\; \;  \rm{as} \; \; n \to \infty.
\ee
\end{theorem}
The main tool is to embed the two-dimensional random walk into a vertical SRW and an horizontal {\em random walk in random scenery} \cite{KS}. The fluctuations of the latter being of order of $n^{3/4}$, this explains why combined with the vertical SRW --whose fluctuations are of order of $n^{1/2}$-- it requires a proper normalization of the order of $n^{5/4}$, see also \cite{CP2,GPLN1,GPLN2}. This strong estimate (\ref{LLT}) implies the convergence of the {\em annealed Green function}
\be\label{UCP}
U_{\rm{cp}}:=\sum_{n=0}^\infty \pee[M_n=(0,0)] < \infty
\ee
which in turns implies this of the  {\em quenched Green function} for $\rho$-a.e. orientation $\eps$ :
\be \label{GreenQue}
0 < U_{\rm{cp}}^{(\eps)}:=\sum_{n=0}^\infty \pee^{(\eps)} [M_n=(0,0)] < \infty,\;  {\rm with} \; 
U_{\rm{cp}}=\E_\rho \big[U_{\rm{cp}}^{(\eps)}\big]>0.
\ee
This also implies\footnote{Although the transience under  in this  quenched law has been proven before, using slightly weaker estimations, but following similar techniques in the vein of Fourier's analysis, see \cite{CP2,GPLN1}.} by Borel-Cantelli the transience of the SRW on $\L^\eps$ for $\rho$-a.e. orientation $\eps$. Thus, the usual dichotomy on $\Z^d$ (P\'olya, 1923) between low dimensions (recurrence for $d=1,2$) and higher dimensions (transience for $d \geq 3$) is broken by the extra-randomness of the orientations\footnote{While it is also proved in \cite{CP2} that deterministic alternate horizontal orientations do not break this recurrence.}. In order to precise the characteristics of this two-dimensional transient random walk, we focus in this paper on the asymptotic behavior of its {\em range} $R_n$, defined to be the number of distinct sites visited by the walker during the first $n$ steps:
$$
R_n={\rm Card} \big\{ M_0, M_1, \dots, M_{n-1} \big\}.
$$
It has been first studied for SRW on $\Z^d$ by Dvoretsky and Erd\"os (\cite{DE}, 1951) who provided estimates of its expectation together with (weak and strong) laws of large numbers under different forms for dimensions $d=2,3,4,\dots$\footnote{Later on, Jain {\em et al}. (\cite{JP2,JP5}, 1970's) established a {\em Central Limit Theorem} (CLT), see Section 5.}.

\section{Results}
\begin{theorem}\label{LinGrowth} The expectation of the range grows linearly :
\beq \label{Qexp}
{\rm For} \; \rho{\rm -a.e.}(\eps), \;     \E^{(\eps)}[R_n] \;  = \; n \cdot\gamma_{\rm{cp}}^{(\eps)} + \circ\big(n\big)
 \; \; {\rm with} \; \; \gamma_{\rm{cp}}^{(\eps)}=(U_{\rm{cp}}^{(\eps)})^{-1} \in \; ]0,1]\\
  \label{Annexp}
\E[R_n] \; = \; n \cdot\gamma_{\rm{cp}} + \circ\big(n\big)
 \; \; {\rm with} \; \; \gamma_{\rm{cp}} =\E_\rho \Big[\frac{1}{U_{\rm{cp}}^{(\eps)}} \Big] \in \; ]0,1].
\eeq
\end{theorem}
The rates of growth $\gamma_{\rm{cp}}$  and $\gamma_{\rm{cp}}^{(\eps)}$ are well-defined as the {\em probability of escape}\footnote{They are related to the notion of capacity of a set reduced to a single point, see \cite{spi3}. The notation $\gamma_{\rm cp}$ stems for Campanino and P\'etritis who first introduced this peculiar random walk in \cite{CP2}.} in next section. 
We emphasize that $\gamma_{\rm cp}$ is {\em not} given by the inverse of the annealed Green function $U_{\rm cp}$, which coincides with the expectation of the quenched Green function $U_{\rm{cp}}^{(\eps)}$. It indeed coincides with {\em the expectation of the inverse of the quenched Green function}  and when the orientations $\eps$ are truly random, these two quantities are not necessarily equal\footnote{This phenomenon occurs rather often in disordered systems or for random walks in random environment.}.\\

One gets thus a linear growth of the expectations of the range similar to the behavior in the space described in \cite{DE}, where a rate $\gamma_3 >0$ is defined similarly, but on a 2-dimensional manifold instead of a 3-dimensional one. The walker visits thus a strictly positive fraction of $n$ sites, on the contrary to the standard planar SRW, for whom the walker typically visits a fraction $\frac{\pi}{\log{n}}$ of $n$ sites,  that goes to zero as $n$ goes to infinity, see \cite{BCR,DE,JP1,LG}. This can be explained by the larger fluctuations, that make the walker escaping from the ball of radius $\sqrt{n}$, and visiting on the way less points already visited. In dimension two, the estimate (2.20) of \cite{DE} yields $\lim_n \frac{\E[R_n]}{n}=0$ but also  the convergence in probability. Here, we also get :

\begin{theorem}\label{WLLN} [Weak Law of Large Numbers (WLLN)] :
\be \label{ALLN}
 \frac{R_n}{n} \stackrel{\mathbb{P}}{\longrightarrow}_n \; \gamma_{\rm{cp}} =\E_\rho \Big[\frac{1}{U_{\rm{cp}}^{(\eps)}} \Big]> 0.
\ee
\end{theorem}

\section{Linear growth of the expected range}
To prove Theorem \ref{LinGrowth}, we follow the road of the original study of \cite{DE}, generalized afterwards by Spitzer \cite{spi3}, and  write this range as a sum of (dependent) random variables $R_n=\sum_{k=0}^{n-1} \mathbf{1}_{A_k}$ where $A_k$ is the event that the walker discovers a new site at the $k^{\rm{th}}$ step i.e.
$$
A_0=\Omega,\; A_k:=\{M_k \neq M_j, \; \forall j=0, \dots, k-1\}.
$$
We also introduce the {\em probability of escape at time k} to be $\gamma_{\rm{cp}}(k):=\pee(A_k)$. As in \cite{DE,spi3}, but with a different manner, we prove that it in fact coincides with the probability that the walk does not come back to its origin during the first $k$ steps.
\bl \label{gammacpn0}
 Denote, for $k \geq 1, \; B_k:=\{M_l \neq (0,0), \; \forall l=1, \dots, k \}$. 
Then $\gamma_{\rm{cp}}(k)=\pee  \big(B_k \big).$
\el

\bpr
On the contrary to the SRW on $\Z^d$, we cannot write $M_n$ as a sum of i.i.d. random variables, but in fact the result can be deduced from the reversibility of the walk. Write $$\pee(A_k)= \sum_{x \in \Z^2} \pee(A_k \cap \{M_k=x\})= \sum_{x \in \Z^2}\E_\rho\big[\pee^{(\eps)}(A_k \cap \{M_k=x\}) \big]$$
and use that for a fixed $\eps$, it corresponds to any trajectory in $A_k$ starting from the origin a unique reversed trajectory in $B_k$, of equal length and equal weight, that is at the origin at $k$: 
\begin{eqnarray*}
 \pee(A_k \cap \{M_k=x\})&=& \pee \big[ \cap_{j=0}^{k-1} \{M_k \neq M_j\} \cap \{M_k=x\} \big]\big]\\
&=&  \sum_{m_j \neq m_k \in \Z^2, j <k}   \E_\rho\big[ \pee^{(\eps)}\big[(M_0, \dots, M_j, \dots, M_k)=(0, \dots, m_j, \dots x)\big]\big]\\
&=& \sum_{m_l \neq m_k \in \Z^2, l <k} \E_\rho \big[\pee^{(-\eps)} \big[(M_0, \dots, M_l, \dots M_k)=(x, \dots, m_l, \dots, 0) \big] \big]\\
&=& \sum_{m_l \neq m_k \in \Z^2, l <k} \E_\rho \big[\pee^{(-\eps)} \big[(M_0, \dots, M_l, \dots M_k)=(0, \dots, m_l, \dots, -x) \big] \big]\\
&=& \pee (B_k \cap \{M_k =-x\})
\end{eqnarray*}
where we use in the last lines the translation-invariance of $\rho$. Integrating out over all the possible final points, one gets
$\pee(A_k)=\sum_{x \in \Z^2} \pee (B_k \cap \{M_k =-x\}) = \pee(B_k)$.
\epr

Hence, the escape probability at time $k$ coincides with the probability of no return to the origin until time $k$. These events $B_k$ are, on the contrary to the $A_k$'s, decreasing events ($B_{k+1} \subset B_k$), in such a way that we get a decreasing sequence $1=\gamma_{\rm{cp}}(1) \geq \dots  \geq \gamma_{\rm{cp}}(k) \geq \gamma_{\rm{cp}}(k+1) \geq \dots \geq 0$. Together with the transience of the walk, this proves that the so-called {\em probability of escape} $\gamma_{\rm{cp}}$ exists and is strictly positive:
$0 < \gamma_{\rm{cp}}:= \lim_k \gamma_{\rm{cp}}(k) \leq \gamma_{\rm{cp}}(k)$ for all $k \geq 0$. We use now the LLT (\ref{LLT}) to get an estimation the growth of the average range, 
\be \label{range}
\E \big[R_n\big]=\sum_{k=0}^{n-1} \pee[A_k]=\sum_{k=0}^{n-1} \gamma_{\rm{cp}}(k).
\ee
Like in \cite{DE}, we partition the paths according to the last return to the origin occurring (strictly) before some given time $n$. The origin can only be reached at even times, so we consider $m=(n-1)/2$ for $n$  even (and $m=n/2 -1$ for $n$ odd) to write, for a given orientation $\eps$,
\be \label{pathsone3}
 \sum_{k=0}^{m}  \pee^{(\eps)} \big[ M_{2k}=(0,0),  M_j \neq (0,0), \; \forall j,\; 2k < j \leq n-1 \big] = 1
\ee
where, by the Markov property of the quenched measure, the summands of (\ref{pathsone3}) are
$$
\pee^{(\eps)} \big[ M_{2k}=(0,0)\big]\cdot \pee^{(\eps)} \big[ M_j\neq(0,0), \; \forall j=2k+1, \dots, n-1 \; \mid M_{2k}=(0,0) \big].
$$ 
Introduce now the following characteristics for the quenched law, for a given orientation $\eps$:
$$
u_{\rm{cp}}^{(\eps)}(k) := \pee^{(\eps)}[M_k=(0,0)] \; {\rm and} \; \gamma_{\rm{cp}}^{(\eps)}(k) := \pee^{(\eps)}[B_k]=\pee^{(\eps)}[M_j\neq (0,0), \forall j, \;  1 <j \leq n].
$$
For $\rho$-a.e. $\eps$, the quenched escape probability $\gamma_{\rm{cp}}^{(\eps)}:=\lim_k \gamma_{\rm{cp}}^{(\eps)}(k) >0$ exists and by symmetry, the probability of discovering a new point at time $k$ is also $\pee^{(\eps)}[A_k]=\gamma_{\rm{cp}}^{(-\eps)}(k)=\gamma_{\rm{cp}}^{(\eps)}(k)$.\\

The techniques developed by \cite{DE} relies on the LLT, here valid in the annealed set-up, yielding the existence of a strictly positive and finite {\em annealed Green function} (\ref{UCP}) and, for $\rho$-a.e$(\eps)$, of a {\em quenched Green function} (\ref{GreenQue})  in such a way that $U_{\rm{cp}}=\E_\rho[U_{\rm{cp}}^{(\eps)}]$. The renewal structure inherited from the Markov property is enough to get 
$$
 \pee^{(\eps)} \big[ M_j\neq(0,0), \; \forall j=2k+1, \dots, n-1 \; \mid M_{2k}=(0,0) \big]= \gamma^{(\eps)}_{\rm{cp}}(n-2k)
$$
so that (\ref{pathsone3}) becomes here, for $\rho$-almost every orientation $\eps$ and for all $n \in \N$
\be\label{pathsone4}
\sum_{k=0}^m u^{(\eps)}_{\rm{cp}}(2k) . \gamma^{(\eps)}_{\rm{cp}}(n-2k)=1
\ee
with $m=(n-1)/2$ for $n$ odd and $m=n/2-1$ for $n$ even. This implies the following

\begin{lemma}\label{LEMAPROUVER}
\begin{enumerate}
\item 
${\rm For} \; \rho{\rm -a.e} \; \eps,\; \gamma_{\rm{cp}}^{(\eps)}.{U_{\rm{cp}}^{(\eps)}}=1 \; {\rm and} \;  \gamma_{\rm{cp}}= \E_\rho\Big[ \frac{1}{U_{\rm{cp}}^{(\eps)}}\Big] >0.$
\item
For all $n \in \N$, there exists $B(n)=\circ(1)$  such that
\be\label{growth3}
0 < \gamma_{\rm{cp}} \leq  \gamma_{\rm{cp}}(n) \leq \gamma_{\rm{cp}} + B(n).
\ee
\end{enumerate}
\end{lemma}
\bpr
Let $\eps$ such that (\ref{GreenQue}) is true, fix $1<l<m$ and split the lhs of (\ref{pathsone4}) to write it
$$\sum_{k=0}^l u^{(\eps)}_{\rm{cp}}(2k). \gamma^{(\eps)}_{\rm{cp}}(n-2k) + \sum_{k=l+1}^m u^{(\eps)}_{\rm{cp}}(2k) .\gamma^{(\eps)}_{\rm{cp}}(n-2k)=1.$$
Use the monotonicity of  $\gamma_{\rm{cp}}^{(\eps)}(k)$ to get a lower bound of the first term of the lhs:
$$
\sum_{k=0}^l u^{(\eps)}_{\rm{cp}}(2k) . \gamma^{(\eps)}_{\rm{cp}}(n-2k) \leq \gamma^{(\eps)}_{\rm{cp}}(n-2l)\cdot \sum_{k=0}^l u^{(\eps)}_{\rm{cp}}(k) 
$$
and the fact that these escape probabilities are indeed probabilities for the second term:
$$
\sum_{k=l+1}^m u^{(\eps)}_{\rm{cp}}(2k) .\gamma^{(\eps)}_{\rm{cp}}(n-2k) \leq  \sum_{k=l+1}^m u^{(\eps)}_{\rm{cp}}(2k) 
$$
to eventually  get the lower bound $\gamma^{(\eps)}_{\rm{cp}}(n-2l) . \sum_{k=0}^l u^{(\eps)}_{\rm{cp}}(2k) \geq 1 - \sum_{k=l+1}^m u^{(\eps)}_{\rm{cp}}(2k).$
Consider now $l\longrightarrow \infty$ such that $n-2l \longrightarrow \infty$  as $n$ goes to infinity, to get  that for $\rho{\rm{-a.e.}} \;  \eps$
$$
\gamma^{(\eps)}_{\rm{cp}}.  U_{\rm{cp}}^{(\eps)} \geq 1
$$
or, the quenched Green function being strictly positive, $\gamma^{(\eps)}_{\rm{cp}}  \geq \frac{1}{U_{\rm{cp}}^{(\eps)}}, \; \rho$-a.s.
By monotonicity one gets in particular for all $n \in \N$ and for $\rho$-a.e. $\eps$
\be \label{controlescquen}
\gamma^{(\eps)}_{\rm{cp}} (n) \geq \frac{1}{U_{\rm{cp}}^{(\eps)}}.
\ee

To get the lower bound, we proceed like in \cite{DE} with a weaker result\footnote{because we do not know whether the quenched local limit theorem is valid or not.} and substract $\frac{1}{U_{\rm{cp}}^{(\eps)}} \sum_{k=0}^{m} u^{(\eps)}_{\rm{cp}}(2k)$ to both sides of (\ref{pathsone4}) to  get first that for $\rho$-a.e. orientation $\eps$,
$$
u_{\rm{cp}}^{\eps}(0)  \cdot \Big(\gamma^{(\eps)}_{\rm{cp}} (n) - \frac{1}{U_{\rm{cp}}^{(\eps)}}\Big) +   \sum_{k=1}^{m} u^{(\eps)}_{\rm{cp}}(2k) . \Big(\gamma^{(\eps)}_{\rm{cp}} (n-2k) - \frac{1}{U_{\rm{cp}}^{(\eps)}}\Big) =  1-\frac{1}{U_{\rm{cp}}^{(\eps)}} \sum_{k=0}^m u^{(\eps)}_{\rm{cp}}(2k)
$$
$$
{\rm so \; that} \; \;\;\;\; \;\;\; u_{\rm{cp}}^{(\eps)}(0)  \cdot \Big(\gamma^{(\eps)}_{\rm{cp}} (n) - \frac{1}{U_{\rm{cp}}^{(\eps)}}\Big)  \; \leq \; 1-\frac{1}{U_{\rm{cp}}^{(\eps)}} \sum_{k=0}^m u^{(\eps)}_{\rm{cp}}(2k).$$
Using (\ref{controlescquen}) and $u_{\rm{cp}}^{(\eps)}(0)=1$, let $n$ (and $m$) going to infinity to get for $\rho$-a.e. $\eps$
 $$
\gamma^{(\eps)}_{\rm{cp}}  \leq \frac{1}{U_{\rm{cp}}^{(\eps)}} \; \;  \; {\rm and  \; thus}  \; \; \gamma^{(\eps)}_{{\rm cp}}  = \frac{1}{U_{{\rm cp}}^{(\eps)}}, \; {\rm and} \; \gamma_{{\rm cp}}  = \E_\rho \Big[ \frac{1}{U_{{\rm cp}}} \Big].
$$
Eventually, we also get that $\rho$-a.s., for all $n \in \N$
$$
0 < \gamma^{(\eps)}_{\rm{cp}} \leq  \gamma^{(\eps)}_{\rm{cp}} (n) \leq  \gamma^{(\eps)}_{\rm{cp}} + B^{(\eps)}(n)
$$
where $B^{(\eps)}(n)= 1- \frac{1}{U_{\rm{cp}}^{(\eps)}} \sum_{k=0}^m u^{(\eps)}_{\rm{cp}}(2k)  = \frac{U_{\rm{cp}}^{(\eps)} - \sum_{k=0}^m u^{(\eps)}_{\rm{cp}}(2k)}{U_{\rm{cp}}^{(\eps)}}$ goes $\rho$-a.s. to $0$.  Taking the expectations w.r.t. $\rho$, this yields the annealed result (\ref{growth3})
where, by dominated convergence, 
$$
B(n)=\E_\rho \Big[\frac{U_{\rm{cp}}^{(\eps)} - \sum_{k=0}^{m} u^{(\eps)}_{\rm{cp}}(2k)}{U_{\rm{cp}}^{(\eps)}} \Big]=\E_\rho\Big[\frac{1}{U_{\rm{cp}}^{(\eps)}}. \sum_{k=m+1}^\infty u_{\rm{cp}}^{(\eps)}(2k)\Big] \; \longrightarrow_n \; 0.
$$
\epr

This provides an estimation of the expected range using  (\ref{range}) to get 
$$n \cdot \gamma_{\rm{cp}} \leq \E[R_n] \leq n \cdot \gamma_{\rm{cp}} + G(n)$$
where by Cesaro's theorem, 
$$
G(n)=\sum_{k=0}^{n-1} B(k)= \sum_{k=0}^{n-1} \E_\rho\Big[\frac{1}{U_{\rm{cp}}^{(\eps)}}. \sum_{l=m(k)+1}^\infty u_{\rm{cp}}^{(\eps)}(2l)\Big]=\circ \big(n\big).
$$
One can proceeds similarly in the quenched case and eventually gets  Theorem \ref{LinGrowth}.

\section{Weak Law of large numbers}

Theorem \ref{LinGrowth} provides thus a linear growth of the expectation of the range 
$$
\frac{\E[R_n]}{n} \; \longrightarrow_n \; \gamma_{\rm{cp}} = \E_\rho \Big[\frac{1}{U_{\rm{cp}}^{(\eps)}} \Big] > 0
$$
similar to the spatial behavior described in \cite{DE} where the limit $\gamma_3 >0$ is defined similarly. This walker goes further than the usual planar one, visiting much more sites but less often. For the standard SRW on the standard (unoriented) version of $\Z^2$, the estimate (2.20) of \cite{DE} 
$$
\E[R_n]= n \cdot \frac{\pi}{\log{n}} + \mathcal{O} \Big(\frac{n \log{\log{n}}}{\log^2{n}}    \Big)
$$
yields $\lim_n \frac{\E[R_n]}{n} = 0$ while  Spitzer \cite{spi3} also proved that $\frac{R_n}{n} \; \stackrel{\mathbb{P}}{\longrightarrow_n} \; 0.$
To investigate this weak LLN\footnote{Established for all $d \geq 2$ in \cite{DE}, who also derive strong LLN.}, we need to estimate the variance of $R_n$, defined to be
\be \label{annvar1}
V_{{\rm cp}}(n):= \sigma^2(R_n) = \E\Big[\big(R_n-\E[R_n]\big)^2 \Big] = \E[R_n^2] - \big(\E[R_n]\big)^2
\ee
which is also the $\rho$-expectation of the quenched variance, defined for a given orientation $\eps$ by
\be \label{annvar2}
V_{{\rm cp}}^{(\eps)}(n):=\E^{(\eps)}\Big[\big(R_n-\E[R_n]\big)^2 \Big] =\E^{(\eps)}[R_n^2] - \big(\E^{(\eps)}[R_n]\big)^2.
\ee
Introduce for all $j<k$ the events $A_{j,k}$ defined as
$$
A_{0,k}=A_k,\; A_{j,k}= \big\{M_k \neq M_l, \forall l=j,\dots, k-1\}.
$$
Re-write now (\ref{annvar1}) and (\ref{annvar2}) as follows
\begin{eqnarray*}
V_{{\rm cp}}(n)&=&\E\big[R_n^2\big] - \Big(\E\big[R_n\big]\Big)^2
= \E\Big[\big(\sum_{j=0}^{n-1} \mathbf{1}_{A_j} \big)^2\Big] - \Big( \E \big[\sum_{j=0}^{n-1}\mathbf{1}_{A_j} \big] \Big)^2\\
&=& \sum_{j,k=0}^{n-1} \Big(\pee \big[ A_j \cap A_k \big]-  \pee\big[A_j\big] . \pee \big[A_{k} \big] \Big)\\
V_{{\rm cp}}^{(\eps)}(n) &=& \sum_{j,k=0}^{n-1} \Big(\pee^{(\eps)} \big[ A_j \cap A_k \big]-  \pee^{(\eps)}\big[A_j\big] . \pee^{(\eps)} \big[A_{k} \big] \Big).
\end{eqnarray*}
Following carefully again the road of \cite{DE} or \cite{spi3}, we establish now the following bound, not optimal\footnote{Investigations around a quenched LLT should lead to $V_{{\rm cp}}(n)=\mathcal{O} \big( n^{3/2} \big)$, see Section 5.} but sufficient to get afterwards a weak law of large numbers:
\bp 
The variance of the range of the SRW on the oriented lattices satisfies
\be\label{varn2}
V_{{\rm cp}}(n)=\circ \big( n^{2} \big) .
\ee
\ep
\bpr
The main ingredient is a sub-additivity of the summands of the variance, that we cannot get using the standard methods of \cite{DE,spi3}. Hence, we first work on the quenched law:
\begin{lemma} \label{subadd}
For all $0\leq j <k$, for all $\eps$, 
\be \label{insubadd}
\pee^{(\eps)}\big[A_j \cap A_k \big] \leq \pee^{(\eps)}\big[A_j\big] . \pee^{(\eps)} \big[ A_{j,k} \big].
\ee
\end{lemma}
\bpr
Use that the quenched law $\pee^{(\eps)}$ is Markov for any orientation $\eps$ to  get for $0 \leq j<k$
\begin{eqnarray*} 
\pee^{(\eps)} [A_j \cap A_k] &=& \pee^{(\eps)} \big[\{M_j \neq M_i, \forall i<j \} \cap \{M_k \neq M_l, \forall l < k \} \big]\\
&\leq& \pee^{(\eps)} \big[\{M_j \neq M_i, \forall i<j \} \cap \{M_k \neq M_l, \forall j \leq l < k \} \big]\\
&=& \pee^{(\eps)} [A_j ] . \pee^{(\eps)} \big[ \{M_k \neq M_l, \forall j \leq l < k \} \big]=  \pee^{(\eps)}\big[A_j\big] . \pee^{(\eps)} \big[ A_{j,k} \big].
\end{eqnarray*}
\epr
\begin{remark}
Inequality (\ref{insubadd}) relies on the Markovian character of the quenched law, not true in the annealed case. Indeed, taking the expectation under $\rho$ in both sides  yields
$$
\pee \big[A_j \cap A_k \big] \leq \E_\rho \Big[ \pee^{(\eps)} \big[A_j \big] . 
\pee^{(\eps)} \big[  A_{j,k} \big] \Big]
$$
and it is an open question whether the product structure of $\rho$ allows to get
$$
 \E_\rho \Big[ \pee^{(\eps)} \big[A_j \big] . \pee^{(\eps)} \big[  A_{j,k}\big] \Big] \leq \E_\rho \Big[ \pee^{(\eps)} \big[A_j \big]\Big] \cdot \E_\rho  \Big[\pee^{(\eps)} \big[  A_{j,k}\big]\Big].
$$

One would get, by translation-invariance of $\rho$, the standard inequality \cite{DE,spi3} because
\be \label{HomAnn}
 \E_\rho \Big[\pee^{(\eps)} \big[  A_{j,k}\big] \Big]= \E_\rho \Big[ \pee^{(\eps)} \big[ A_{k-j} \big]\Big]=\pee \big[A_{k-j} \big].
\ee
 
\end{remark}
Using now the estimate (\ref{insubadd}) and the expression (\ref{range}), we can estimate (\ref{annvar2})
\begin{eqnarray*}
V_{{\rm cp}}^{(\eps)}(n)&=&2 \sum_{j=0}^{n-1}\sum_{k=j+1}^{n-1} \Big(\pee^{(\eps)} \big[ A_j \cap A_k \big]-  \pee^{(\eps)}\big[A_j\big] . \pee^{(\eps)} \big[ A_{k} \big]\Big) + \sum_{j=0}^{n-1} \big(\pee^{(\eps)} \big[ A_j \big]-\pee^{(\eps)} \big[ A_j \big]^2 \big)\\
&\leq & 2 \sum_{j=0}^{n-1} \pee^{(\eps)} \big[ A_j \big] \cdot \sum_{k=j+1}^{n-1}  \Big(\pee^{(\eps)} \big[ A_{j,k} \big] - \pee^{(\eps)} \big[ A_k \big] \Big) + \sum_{j=0}^{n-1} \pee^{(\eps)} \big[ A_j \big]
\end{eqnarray*}
so that
$$
\frac{1}{n^2} V_{{\rm cp}}^{(\eps)}(n) \leq  \frac{2}{n} \sum_{j=0}^{n-1} \Big( \pee^{(\eps)} \big[ A_j \big] \cdot \sum_{k=j+1}^{n-1}  \frac{1}{n}\Big(\pee^{(\eps)} \big[ A_{j,k} \big] - \pee^{(\eps)} \big[ A_k \big] \Big) \Big) + \frac{\E^{(\eps)} \big[R_n \big]}{n^2} = G_n(\eps) + \frac{\E^{(\eps)} \big[R_n \big]}{n^2}.
$$
The last term of the rhs goes $\rho$-a.s. to zero by (\ref{Qexp}) while we write
$$
G_n(\eps)=2 \gamma_{\rm cp} \cdot \E_\rho \Big[\frac{1}{n} \sum_{j=0}^{n-1}  \sum_{k=j+1}^{n-1} \frac{1}{n}  \Big(\pee^{(\eps)} \big[ A_{j,k} \big] - \pee^{(\eps)} \big[ A_k \big] \Big) \Big)\Big] + 2 D_n(\eps)=2 (F_n(\eps) + D_n(\eps))
$$ 
in such a way that we control the annealed variance by 
$\frac{1}{n^2} V_{{\rm cp}}(n) \leq 2\E_\rho[F_n] + 2\E_\rho[D_n]$.\\

To deal with the second term, remark that, for given $j$ and $k$,  $A_{j,k}= A_k \cup \tilde{A}_{j,k}$ where the events $\tilde{A}_{j,k}$ consists of the trajectories visiting at $k$ a point not visited since $j$ but who has been visited before. In particular, since $A_{j,k} \subset A_{0,k}=A_k$, 

$$
0 \leq \frac{1}{n} \sum_{k=j+1}^{n-1}  \Big(\pee^{(\eps)} \big[ A_{j,k} \big] - \pee^{(\eps)} \big[ A_k \big] \Big) = \frac{1}{n} \sum_{k=j+1}^{n-1} \pee^{(\eps)} \big[ \tilde{A}_{j,k} \big]\leq 1
$$
so that
$$
0 \leq \E_\rho[D_n] \leq \E_\rho \Big[ \frac{1}{n} \sum_{j=0}^{n-1} \Big( \pee^{(\eps)} \big[ A_j \big] - \gamma_{\rm{cp}} \Big)\Big] = \E_\rho \Big[ \Big(\frac{1}{n} \sum_{j=0}^{n-1} \pee^{(\eps)} \big[ A_j \big] \Big) - \gamma_{\rm{cp}} \Big] = \frac{\E[R_n]}{n} - \gamma_{\rm{cp}}
$$
that goes to zero by (\ref{Annexp}). To deal with $F_n$, we write
$$
\E_\rho[F_n] \leq \gamma_{\rm{cp}} \cdot \Big[ \sum_{j=0}^{n-1} \frac{1}{n} \sum_{k=j+1}^{n-1}  \Big(\pee \big[ A_{j,k} \big] - \pee \big[ A_k \big] \Big) \Big]\leq  \sum_{j=0}^{n-1} \frac{1}{n} \gamma_{\rm{cp}} \cdot \max_{j=0,\dots k-1} \sum_{k=j+1}^{n-1}  \Big(\pee \big[ A_{k-j} \big] - \pee \big[ A_k \big]\Big)
$$
where to get the last inequality we have used (\ref{HomAnn}) for the annealed measure\footnote{This step is not true in the quenched case so we cannot get the same bound, at least in this way.}. Now, we can work exactly like in the standard case treated in \cite{DE,spi3}: The balance between the number of possible points to discover and the number of points already visited reaches its maximum for $j=\big[\frac{n}{2}\big]$ so that,  using (\ref{Annexp}), we get (\ref{varn2}) because
$$
0 \leq \E_\rho[F_n] \leq \gamma_{\rm cp} \cdot \frac{1}{n} \cdot \Big (\E \Big[ R_{n-[n/2]}+ R_{[n/2]} - R_n \Big]\Big) \leq \gamma_{\rm cp} \cdot \Big(\frac{1}{2} \gamma_{\rm{cp}} + \frac{1}{2} \gamma_{\rm{cp}} - \gamma_{\rm{cp}} \Big) +  \circ (1 )= \circ (1 ).
$$
 \epr

Using Markov's inequality one gets Theorem \ref{ALLN}, because for all $\delta >0$
 $$
\pee \Big[ \big|\frac{R_n}{n} - \gamma_{\rm{cp}} \big| > \delta \Big] \leq  \frac{1}{n^2 \delta^2} .\E\big[ |R_n - n . \gamma_{\rm{cp}} | \big] 
\leq \frac{1}{n^2 \delta^2} . V_{{\rm cp}}(n) + \frac{1}{\delta^2} . \Big(\gamma_{\rm{cp}} - \E \Big[\frac{R_n}{n}\Big] \Big)^2
$$
and the WLLN in the annealed set-up, by (\ref{varn2}) and (\ref{growth3}). As a by-product, one recovers also in the quenched WLLN for $\rho$-a.e. orientation.

\section{Conclusions and perspectives}

Further investigations, in the spirit of Jain {\em et al.} \cite{JP2,JP5}, would require a quenched local limit theorem or at least more accurate asymptotic of the variance of the range, using e.g. a less crude inequality than (\ref{insubadd}), and in this transient case the relationship between the range and the number of points that are never revisited.  We suspect that in fact the variance is of order $n^{3/2}=n \sqrt{n} $, and that this should lead to an unconventional CLT:
$$
\frac{R_n -n \gamma_{\rm{cp}}}{\sqrt{n \sqrt{n}}} \; \stackrel{\mathcal{L}}{\Longrightarrow} \; \mathcal{N}(0,1)
$$
like in the three-dimensional case (where the normalization is $\sqrt{n \ln{n}})$, while in the two-dimensional case  the limiting law is the so-called  self-intersection local times \cite{LG}.\\

{\bf Aknowledegments :} I am grateful to Jean-Baptiste Bardet (Rouen), Frank den Hollander (Leiden) and Bruno Schapira (Orsay) for their interest and their advices.


\end{document}